\documentclass{gtart_h}

\def\ifplaintex{\expandafter\ifx\csname documentclass\endcsname\relax}

\def\gtp{{\mathsurround=0pt\it $\cal G\mskip-2mu$eometry \&\ 
$\cal T\!\!$opology $\cal P\!$ublications}}  

\def\recd{{\small Received:\qua\receiveddate\ifx\reviseddate\relax
\else\qquad Revised:\qua\reviseddate\fi\par}} 

\def\lognumber#1{\def\thelognumber{#1}}
\def\volumenumber#1{\def\thevolumenumber{#1}}
\def\volumeyear#1{\def\thevolumeyear{#1}}
\def\papernumber#1{\def\thepapernumber{#1}}
\def\pagenumbers#1#2{\def\startpage{#1}\def\finishpage{#2}}
\def\published#1{\def\publishdate{#1}}

\def\received#1{\def\receiveddate{#1}}
\def\revised#1{\def\reviseddate{#1}}
\def\accepted#1{\def\accepteddate{#1}}

\def\asciiemail#1{\def\theasciiemail{#1}}

\long\def\asciiabstract#1{\long\def\theasciiabstract{#1}}
\def\asciikeywords#1{\def\theasciikeywords{#1}}

\let\\\par\let\thelognumber\relax\let\thevolumenumber\relax
\let\thepapernumber\relax\let\thevolumeyear\relax\let\startpage\relax
\let\finishpage\relax\let\publishdate\relax\let\receiveddate\relax
\let\reviseddate\relax\let\accepteddate\relax\let\theasciititle\relax
\let\theasciiauthors\relax
\let\theasciiabstract\relax\let\theasciikeywords\relax

\let\theasciiemail\relax

\ifplaintex
\font\logobig=cmssbx10 scaled 3836
\font\logomed=cmssbx10 scaled 2557
\else
\font\logobig=cmssbx10 scaled 4200
\font\logomed=cmssbx10 scaled 2800
\fi

\long\def\makeagttitle{   
\count0=\startpage
\agt\hfill      
\hbox to 45truept{\vbox to 0pt{\vglue -13truept{\logomed A\kern -.37em{\logobig 
T}\kern -.38em G}\vss}\hss}
\break
{\small Volume \thevolumenumber\ (\thevolumeyear)
\startpage--\finishpage\nl
Published: \publishdate}

\vglue .25truein

{\parskip=0pt\leftskip 0pt plus
1fil\def\\{\par\smallskip}{\Large\bf\thetitle}\par\medskip} \vglue
0.05truein

{\parskip=0pt\leftskip 0pt plus 1fil\def\\{\par}{\sc\theauthors}
\par\medskip}%
 
\vglue 0.03truein 

{\small\leftskip 25truept\rightskip 25truept{\bf Abstract}\stdspace\theabstract

{\bf AMS Classification}\stdspace\theprimaryclass
\ifx\thesecondaryclass\relax\else; \thesecondaryclass\fi\par
{\bf Keywords}\stdspace \thekeywords\par}\vglue 7truept

}   

\ifplaintex
\font\phead=cmsl9 scaled 950
\font\pnum=cmbx10 scaled 913
\font\pfoot=cmsl9 scaled 950
\headline{\vbox to 0pt{\vskip -4.5mm\line{\small\phead\ifnum
\count0=\startpage ISSN 1472-2739 (on-line) 1472-2747 (printed)
\hfill {\pnum\folio}\else\ifodd\count0\def\\{ }%
\ifx\theshorttitle\relax\thetitle\else\theshorttitle\fi\hfill{\pnum\folio}
\else\def\\{ and }{\pnum\folio}\hfill\ifx\theshortauthors\relax\theauthors
\else\theshortauthors\fi\fi\fi}\vss}}
\footline{\vbox to 0pt{\vglue 0mm\line{\small\pfoot\ifnum\count0=\startpage
\copyright\ \gtp\hfill\else
\agt, Volume \thevolumenumber\ (\thevolumeyear)\hfill\fi}\vss}}
\else
\font=cmsl9 scaled 1050
\font\lnum=cmbx10 
\font=cmsl9 scaled 1050
\makeatletter
\def\@oddhead{{\small\lhead\ifnum\count0=\startpage ISSN 1472-2739 
(on-line) 1472-2747 (printed)\hfill {\lnum\number\count0}\else\ifodd\count0
\def\\{ }\ifx\theshorttitle\relax \thetitle \else\theshorttitle\fi\hfill
{\lnum\number\count0}\else\def\\{ and }{\lnum\number\count0}
\hfill\ifx\theshortauthors\relax 
\theauthors\else\theshortauthors\fi\fi\fi}}\def\@evenhead{\@oddhead}
\def\@oddfoot{\small\lfoot\ifnum\count0=\startpage\copyright\ \gtp\hfill\else
\agt, Volume \thevolumenumber\ (\thevolumeyear)\hfill\fi}
\def\@evenfoot{\@oddfoot}
\makeatother
\fi
\let\maketitlepage\makeagttitle
\let\makeshorttitle\maketitlepage
\let\maketitle\maketitlepage

\newwrite\gtoutfile
\long\gdef\makeheadfile{  
{\def\\{, }\def\s{ }
\immediate\openout\gtoutfile head.xxx
\immediate\write\gtoutfile{Proxy-for: \ifx\theasciiauthors\relax
\theauthors\else\theasciiauthors\fi\s<\ifx\theasciiemail\relax\theemail\else\theasciiemail\fi>}
\immediate\write\gtoutfile{\noexpand\\}
\immediate\write\gtoutfile{Authors: \ifx\theasciiauthors\relax
\theauthors\else\theasciiauthors\fi}
{\def\\{ }\immediate\write\gtoutfile{Title: \ifx\theasciititle\relax
\thetitle\else\theasciititle\fi}}
\immediate\write\gtoutfile{Subj-class: GT or SG, GR etc}
\immediate\write\gtoutfile{MSC-class: \theprimaryclass\ifx\thesecondaryclass\relax\else, \thesecondaryclass\fi}
\immediate\write\gtoutfile{Journal-ref: Algebr. Geom. Topol. \thevolumenumber\s
(\thevolumeyear) \startpage-\finishpage}
\immediate\write\gtoutfile{Comments: Published by Algebraic and
Geometric Topology at}
\immediate\write\gtoutfile{\s\s\s  http://www.maths.warwick.ac.uk/agt/AGTVol\thevolumenumber/agt-\thevolumenumber-\thepapernumber.abs.html}
\immediate\write\gtoutfile{\noexpand\\}
\immediate\write\gtoutfile{}
\ifx\theasciiabstract\relax
\immediate\write\gtoutfile{\theabstract}\else
\immediate\write\gtoutfile{\theasciiabstract}\fi
\immediate\write\gtoutfile{}
\immediate\write\gtoutfile{\noexpand\\}
\immediate\write\gtoutfile{}
\immediate\closeout\gtoutfile}}  

\def\maketitlepage{\makeagttitle\makeheadfile}
\let\makeshorttitle\maketitlepage
\let\maketitle\maketitlepage

\lognumber{35}
\volumenumber{4}
\volumeyear{2004}
\papernumber{35}
\pagenumbers{813}{827}
\received{29 January 2004} 
\revised{4 July 2004}
\accepted{23 September 2004}
\published{23 September 2004}

\usepackage{amsmath,amscd,amssymb}
\usepackage[all]{xy}

\def\sh#1{\subsection*{#1}\addcontentsline{toc}{subsection}{#1}}

\newtheorem{thm}{Theorem}[section]  
\newtheorem*{un-no-thm}{Theorem}
\newtheorem{cor}[thm]{Corollary}     
\newtheorem{lem}[thm]{Lemma}         
\newtheorem{prop}[thm]{Proposition}  

\newtheorem{bigthm}{Theorem}
\newtheorem{bigcor}[bigthm]{Corollary}

\theoremstyle{definition}
\newtheorem{defn}[thm]{Definition}   

\theoremstyle{definition}
\theoremstyle{definition}

\theoremstyle{remark}

\newtheorem*{acks}{Acknowledgements}
\newtheorem*{out}{Outline}

\newtheorem*{rem_no}{Remark}

\newtheorem*{convent}{Conventions}

\def\Top{\bold T\bold o \bold p}

\def\:{\colon\thinspace}

\def\Bbb{\mathbb}
\def\bold{\mathbf}

\def\cal{\mathcal}

\begin{document}
\title{On the homotopy invariance of\\configuration spaces}
\author{Mokhtar Aouina\\John R. Klein}
\address{Department of Mathematics, Wayne State University\\Detroit, 
MI 48202, USA}

\asciiemail{aouina@math.wayne.edu, klein@math.wayne.edu}

\gtemail{\mailto{aouina@math.wayne.edu}\qua
{\rm and}\qua \mailto{klein@math.wayne.edu}}

\begin{abstract} For a closed PL manifold
$M$, we consider the configuration space $F(M,k)$
of ordered $k$-tuples of distinct  points in $M$.
We show that a suitable iterated suspension of $F(M,k)$
is a homotopy invariant of $M$. The number of suspensions
we require depends on three parameters: the number of points $k$,
the dimension of $M$ and the connectivity of $M$.
Our proof uses a mixture
of Poincar\'e embedding theory and fiberwise algebraic topology.
\end{abstract}

\asciiabstract{%
For a closed PL manifold M, we consider the configuration space F(M,k)
of ordered k-tuples of distinct points in M.  We show that a suitable
iterated suspension of F(M,k) is a homotopy invariant of M. The number
of suspensions we require depends on three parameters: the number of
points k, the dimension of M and the connectivity of M.  Our proof
uses a mixture of Poincare embedding theory and fiberwise algebraic
topology.}

\asciikeywords{Configuration space, fiberwise suspension, embedding up to 
homotopy, Poincare embedding}
\maketitle

\section{Introduction}

For a closed PL manifold $M$ and an integer $k \ge 2$, 
we will consider the configuration space
$$
F(M,k) := \{(x_1,...,x_k)|\, x_i \in M \text{ and } x_i \ne x_j 
\text{ for } i \ne j \}\, .
 $$
A fundamental unsolved problem  about these spaces concerns
their {\it homotopy invariance:} when $M$ and
$N$ are homotopy equivalent, is it true that $F(M,k)$ and
$F(N,k)$ are homotopy equivalent?

Here is some background. 
It is known that the based loop space $\Omega F(M,k)$, 
is a homotopy invariant (see Levitt \cite{Levitt}).
When $M$ is smooth, the cohomology of $F(M,k)$ with field
coefficients has been intensively studied 
(see e.g., B\"odigheimer-Cohen-Taylor \cite{BCT}).
When $M$ is a smooth projective variety over $\Bbb C$,
Kriz \cite{Kr} has shown that the rational homotopy type
of $F(M,k)$ depends only on the rational cohomology ring
of $M$.  

When  $k = 2$ we have $F(M,2) = M \times M - \Delta$ is
the deleted product. Even in this instance,
the homotopy invariance question is still not completely settled.
However, shortly after the first draft of this paper was circulated, 
R.\ Longoni and P.\ Salvatore partially settled the
question by showing that the deleted product spaces 
of the homotopy equivalent (but not {\it simple} homotopy
equivalent) lens spaces $L(7,1)$ and $L(7,2)$
have distinct homotopy types (see \cite{LS}).

The purpose of this paper is to show that a suitable iterated suspension
of $F(M,k)$ is a homotopy invariant. The bound
on the number of suspensions we need to
take depends on three parameters: the number of
points, the dimension of $M$
and the connectivity of $M$.

For an unbased space $Y$, we define its $j$-fold (unreduced) suspension 
$$
\Sigma^jY := (*\times S^{j{-}1}) \cup (Y \times D^j)\, ,
$$
 where the union is
amalgamated along $Y \times S^{j{-}1}$ (up to homotopy,
$\Sigma^j Y$ is the {\it join} of $Y$ and $S^{j{-}1}$).

Our main result is

\begin{bigthm}\label{main-thm} Let $M$ and $N$ be 
homotopy equivalent closed {\rm PL} manifolds
of dimension $d$. Assume $M$ is $r$-connected for some $r\ge 0$. 
Then there is a homotopy
equivalence
$$
\Sigma^{\alpha(k,d,r)} F(M,k) \,\, \simeq \,\, \Sigma^{\alpha(k,d,r)} F(N,k) \, ,
$$
where $\alpha(k,d,r) := \max((k-2)d - r + 3,0)$. 
\end{bigthm}

\begin{rem_no} 
Cohen and Taylor (unpublished manuscript)
prove by very different methods that 
the configuration spaces
of {\it smooth} manifolds are stable homotopy invariant.
In their work  the bound on the number suspensions required 
to achieve homotopy invariance is significantly weaker.

Nevertheless, an advantage of their approach is
its applicability to other kinds of configuration spaces.
For example, their results apply as well to 
the {\it unordered}  configuration spaces of a 
smooth manifold. We are unable to analyze 
the latter using our methods.
\end{rem_no}

\begin{bigcor} \label{stable} Let $M$ be a connected
closed {\rm PL} manifold. Then the suspension spectrum  $\Sigma^\infty F(M,k)_+$
is a homotopy invariant of $M$.
\end{bigcor} 

Theorem \ref{main-thm}
can be improved by one dimension
provided that the input manifolds are {\it simple homotopy
equivalent:} 

\begin{bigthm}\label{main-add} With the additional
assumption that $M$ and $N$ 
are simple homotopy equivalent, there
is a homotopy equivalence
$$
\Sigma^{\beta(k,d,r)} F(M,k) \,\, \simeq \,\, \Sigma^{\beta(k,d,r)} F(N,k) \, ,
$$
where $\beta(k,d,r) := \max((k-2)d - r + 2,0)$.
\end{bigthm} 

\begin{rem_no} For simply connected manifolds, a 
 homotopy equivalence is also simple. Thus
Theorem \ref{main-add} improves upon Theorem \ref{main-thm} 
in the $1$-connected case.
\end{rem_no}

The following corollary extends
the work of Levitt \cite{Levitt}, who considered only the case 
of $2$-connected manifolds.

\begin{bigcor} If $M$ is connected, then 
$\Sigma^3 F(M,2)$ is a homotopy
invariant of $M$. Furthermore, if $M$ is $1$-connected, then
$\Sigma F(M,2)$ is a homotopy invariant of $M$.
\end{bigcor} 

(The first part is just a special case of Theorem \ref{main-thm}, whereas 
the second part is a special case of Theorem \ref{main-add}.)

\begin{convent}
We work in the category $\Top$ of compactly generated topological
spaces.

A non-empty space is always $(-1)$-connected.
A space is $0$-connected if it is path connected.
It is $r$-connected for $r >0$ if it is path connected and
its homotopy groups (with respect to a choice
of basepoint) vanish in degrees $\le r$.
 A map $A \to B$
of  spaces is $r$-connected 
if its homotopy fiber at all basepoints
is $(r{-}1)$-connected.
A weak (homotopy) equivalence is an $\infty$-connected map.
If two spaces $A$ and $B$ are related by a chain of weak equivalences,
we will often indicate it by writing $A\simeq B$. A space is {\it homotopy
finite} if it is homotopy equivalent to a finite cell complex.
\end{convent}

\begin{out} In \S2 we describe the construction
of fiberwise suspension and deduce some elementary properties of it.
In \S3 we review the Stallings-Wall theory of embeddings up to homotopy.
\S4 is about decompressing embeddings up to homotopy so as to increase their
codimensions. A key result of this
section concerns the iterated suspension of the complement
of an embedding up to homotopy.  In \S5 Theorem \ref{main-add}
is proved using the Browder-Casson-Sullivan-Wall theorem and
the Stallings-Wall embedding theorem. In \S6 Theorem \ref{main-thm}
is proved using the second author's previous work on Poincar\'e embeddings.
\end{out}

\begin{acks}  We are indebted to the referee for pointing
out a mistake we made when  applying the Browder-Casson-Sullivan-Wall theorem
in a previous version of this paper.
The mistake evaporates when one assumes 
simple homotopy equivalences between the input manifolds:
Theorem \ref{main-add} is an artifact of our original
(erroneous) proof of Theorem \ref{main-thm}.

The proof of Theorem \ref{main-thm} contained here uses a result of the second
author on concordances between Poincar\'e embeddings; the latter 
accounts for the loss of one dimension in the statement of 
Theorem \ref{main-thm}.

The first author is supported by a Wayne State University Rumble Fellowship.
The second author is partially supported by  NSF Grant DMS-0201695.
\end{acks}

\section{Fiberwise suspension}

Let $A\to X$ be a  map of spaces. 
Define $$\Top_{A\to X}$$ to be the category of 
spaces ``between $A$ and $X$.'' Specifically, an {\it object}
  is a space $Y$ and a choice of factorization
$A\to Y \to X$. 
A {\it morphism} is a map of spaces
which is compatible with their given factorizations.
Call a morphism a {\it weak equivalence} if it
is a weak homotopy equivalence of underlying spaces.

We write $\Top_{/X}$ for $\Top_{\emptyset \to X}$.
If $Y \in \Top_{/X}$ is an object, define its {\it (unreduced) $j$-fold
 fiberwise suspension} by
$$
\Sigma^j_X Y \,\, := \,\, (Y \times D^j) \cup (X \times S^{j-1})\, ,
$$
where the union is amalgamated over $Y {\times}S^{j-1}$. 
With respect to the first factor projection map $X \times S^{j{-}1} \to X$,
we get a functor
$$
\Sigma^j_X\: \Top_{/X} \to \Top_{X{\times} S^{j-1} \to X}\, .
$$

\begin{lem} \label{collapse} Let $Y$ and $Z$ be objects of $\Top_{/X}$
whose underlying spaces are path connected
and have the homotopy type of CW complexes. 

Assume for 
some $j \ge 0$ that  $\Sigma^j_XY$ and $\Sigma^j_X Z$ are
weak equivalent objects.
Then there is a weak  equivalence of spaces
$$
\Sigma^j Y \,\, \simeq \,\, \Sigma^j Z\, . 
$$
\end{lem}

\begin{proof} 
The statement is obviously true for $j =0$, so
we will assume that $j > 0$. Moreover,
we may assume that we are given a weak equivalence
$\Sigma^j_X Y \overset\sim\to \Sigma^j_X Z$.

For any object $T\in \Top_{/X}$, we have a cofibration sequence
of spaces
$$
X\times S^{j-1} \to  \Sigma^j_X T \to \Sigma^j (T_+) \, ,
$$ 
where we use that $\Sigma^j (T_+)$ is
$T \times D^j$ with $T\times S^{j{-}1}$ collapsed to
a point. Using this  cofiber sequence for both $Y$ and $Z$, we get a
commutative diagram
$$
\SelectTips{cm}{}
\xymatrix{
\Sigma^j_X Y \ar[r]\ar[d]_{\simeq}  &\Sigma^j (Y_+)\ar[d] \\
\Sigma^j_X Z \ar[r] &\Sigma^j (Z_+)
}
$$
which is also homotopy pushout. It is well-known that
cobase change preserves weak equivalences (see e.g., Hirschhorn 
\cite{Hirschhorn}), 
so it follows that the map $\Sigma^j (Y_+) \to \Sigma^j (Z_+)$ is a weak
 equivalence.

If $A$ is a space, let $A_+$ be the union of $A$ with a disjoint basepoint.
Then one has a weak equivalence 
$$
\Sigma^j (A_+) \, \simeq \, (\Sigma^j A) \vee S^j\, ,
$$
for $j > 0$, where the left side is to be regarded
as the reduced suspension of a based space.

An explicit weak equivalence can be constructed as follows: 
recall that 
$$
\Sigma^j A \,\, = \,\, 
(A \times D^j)\cup_{A\times S^{j-1}} (*\times S^{j-1}) \, .
$$ 
The effect of collapsing  $*\times S^{j-1}$ to a point defines a based map
$$
i\:\Sigma^j A \to \Sigma^j(A_+)\, .
$$ 
Choosing a basepoint for $A$ yields a based map $S^0 \to A_+$.
Let $k\:S^j \to \Sigma^j (A_+)$ 
denote its  $j$-fold suspension. Then we obtain
a map
$$
\begin{CD}
(S^j A) \vee S^j @> i \vee k >> \Sigma^j (A_+)\vee
\Sigma^j (A_+) @>\text{fold} >> \Sigma^j (A_+) \, .
\end{CD}
$$
It is straightforward to check
that this map is a weak equivalence. 

Consequently,
we have weak equivalences
$\Sigma^j (Y_+) \simeq (\Sigma^j Y) \vee S^j$  
and $\Sigma^j (Z_+) \simeq (\Sigma^j Z) \vee S^j$ for $j> 0$.
It follows that
there is a weak equivalence 
$$
 (\Sigma^j Y) \vee S^j \,\, \simeq  \,\,  (\Sigma^j Z) \vee S^j \, .
$$
Because $Y$ and $Z$ are connected, we have that $\Sigma^j Y$ and $\Sigma^j Z$
are $j$-connected. Using Lemma \ref{split} below, we conclude
that the composite
$$
\begin{CD}
\Sigma^j Y @> \text{include} >>  (\Sigma^j Y) \vee S^j 
\simeq   (\Sigma^j Z) \vee S^j 
@>\text{project} >>  (\Sigma^j Z)
\end{CD}
$$ 
is a weak equivalence.
\end{proof}

\begin{lem}\label{split} 
Let $U$ and $V$ be $j$-connected spaces with $j \ge 0$.
Assume $U$ and $V$ are equipped with non-degenerate basepoints.
Assume $h\:{U \vee S^j \to V \vee S^j}$
is a weak  equivalence. Then the composition
$$
\begin{CD}
g\: U @> \text{\rm include} >> U \vee S^j
@> h >> V \vee S^j @> \text{\rm project} >> 
V
\end{CD}
$$
is also a weak equivalence.
\end{lem}

\begin{proof} 
Without loss in generality we can assume that
$U$ and $V$ are CW complexes with no cells in positive
dimensions $\le j$.
By cellular approximation, we may also assume that $h$ is a cellular map.
Then $h$ preserves $j$-skeleta, so there is a
commutative diagram
$$
\SelectTips{cm}{}
\xymatrix{
S^j \ar[r]^{\subset\quad } \ar[d]_{h_{|S^j}} & U \vee S^j \ar[d]^{h}\\
S^j \ar[r]_{\subset\quad } & V \vee S^j \, ,
} 
$$
and it is straightforward to check that 
the left vertical map is a homotopy equivalence.
We infer that the map $U \to V$ obtained by
taking cofibers horizontally is  also a weak equivalence.
But this map coincides with $g$.
\end{proof}

\section{Embeddings up to homotopy}

Let $K$ be a space 
which is homeomorphic to a connected finite complex of dimension $\le k$.
Let $M$ be a PL manifold of dimension $d$, possibly with boundary.
Fix a map $f\: K\to M$.

\begin{defn} An {\it embedding up to homotopy} of $f$ is 
a pair 
$$
(N,h)
$$ in which
\begin{itemize}
\item $N$ denotes a compact codimension zero PL 
submanifold of the interior of $M$;
\item the pair $(N,\partial N)$ is $(n{-}k{-}1)$-connected;
\item $h\:K \to N$ is a simple homotopy equivalence such that
composition $$K \overset h \to N \subset M$$ is homotopic to $f$.
\end{itemize}

A {\it concordance} of embeddings up
to homotopy $(N_0,h_0)$ and $(N_1,h_1)$
of $f$ consists of a locally flat embedding 
$$
e\:N_0 \times [0,1] \subset M \times [0,1]
$$
and a homotopy $H_t\: K \to N_0$
such that 
\begin{itemize}
\item $e$ restricted to $N_0 \times 0$ is the inclusion
and $e$ maps $N_0\times 1$ homeomorphically onto $N_1$. 
\item $H_0 = h_0$ and $H_1$ followed by $e(\cdot,1)$ coincides
with $h_1$.  
\end{itemize}
\end{defn}

\begin{thm}[Stallings \cite{Stallings}, Wall \cite{Wall4}] \label{stallings-wall} Assume 
$\dim K\le k \le d {-}3$. If $f\:K\to M$ is
$(2k {-} d {+} 1)$-connected, then $f$ embeds up to homotopy. Furthermore,
any two embeddings up to homotopy of $f$ are concordant whenever 
$f$ is $(2k {-} d {+} 2)$-connected.
\end{thm}

\section{Decompression}

Let $(N,h)$ be an embedding up to homotopy of $f\: K \to M$.
If $C$ denotes the closure of the complement of $N$ inside $M$, then
$C$ is an object of $\Top_{\partial M \subset M}$.

\begin{defn} The object 
$$
C \in  \Top_{\partial M \subset M}
$$
is called the {\it complement} of $(N,h)$.
\end{defn}

By considering the inclusion $M \times 0 \subset M \times D^j$,
and taking a compact regular neighborhood of $N$ in $M \times D^j$,
we have an associated embedding up to homotopy of the composite
$$
f_j\:K \overset f \to M = M \times 0 \subset M \times D^j\, .
$$
Denote this embedding up to homotopy by $(N_j,h_j)$, where
$N_j \cong N \times D_{1/2}^j$ (here $D_{1/2}^j\subset D^j$ is 
the disk of radius $1/2$)
and $h_j$ is identified with $h$ followed by the inclusion
$N \times 0 \subset N \times D^j_{1/2}$.
 This new embedding up
to homotopy is the {\it $j$-fold decompression} of $(N,h)$.
Note that its complement has the structure of an object of 
$\Top_{\partial(M {\times} D^j)\subset M {\times} D^j}$.

However, to avoid technical problems, we will henceforth
regard the complement as a space over $M$ by projecting away
from the $D^j$ factor.
That is, we will think of the complement
as an object of $\Top_{\partial(M {\times} D^j)\to M}$.

\begin{lem}[Compare {\cite[\S2.3]{Klein_haef}}] \label{decomp}  Assume that $M$ is closed.
Then the complement of
$(N_j,h_j)$  is weak equivalent to
the object
$$
\Sigma^j_M C \,. 
$$
\end{lem}

\begin{proof} The regular neighborhood $N_j$ can be chosen
as $N \times D^j_{1/2} \subset M \times D^j$. 
The complement of $(N_j,h_j)$
is then 
$$
(M \times D^j) - \text{int}(N \times D^j_{1/2})
\,\, =  \,\, 
C \times D^j_{1/2} \,\, \cup \,\, M \times D^j_{[1/2,1]} \, ,
$$
where $D^j_{[1/2,1]}$ denotes the annulus consisting 
of points in $D^j$ whose norm varies
between $1/2$ and $1$. The above union is amalgamated over
$C \times \partial D^j_{1/2}$. 

The subspace of the complement given by
$(C \times D^j_{1/2})\cup (M \times \partial D^j_{1/2})$ 
is evidently isomorphic
to $\Sigma^j_M C$. The inclusion map of this
subspace is, up to isomorphism,
a morphism of  
$\Top_{M {\times} S^{j{-}1}\to M}$.
Furthermore, this inclusion is a weak homotopy equivalence of 
the underlying spaces.
\end{proof}

We conclude this section with a result about the homotopy
type of the iterated suspension of the complements of embeddings up to homotopy. This will be a key ingredient of the proof of Theorem \ref{main-add}.

\begin{prop}\label{susp-comp} Assume
$f\:K^k \to M^d$ is an $r$-connected map, where $M$ is
a closed connected {\rm PL} manifold of dimension $d$, and $\dim K\le  k\le d{-}3$.
Suppose that $f$ has two embeddings up to homotopy 
$(N,h)$ and $(N',h')$ with respective complements
$C$ and $C'$. Then there is a homotopy equivalence, 
$$
\Sigma^j C \,\, \simeq \,\,  \Sigma^j C' \, ,
$$
where $j = \max(2k - d - r + 2,0)$. 
\end{prop}

\begin{proof}
By the Stallings-Wall theorem (\ref{stallings-wall}), with $j = 
 \max(2k - d - r + 2,0)$, we see that the $j$-fold decompressions
of $(N,h)$ and $(N',h')$ are concordant.
Furthermore, it is evident from the 
definitions that concordant embeddings up to homotopy
have homotopy equivalent complements.

Using
Lemma \ref{decomp} we infer that there is a weak equivalence
of objects 
$$
\Sigma^j_M C \,\,  \simeq \,\, \Sigma^j_M C' \, .
$$
By Lemma \ref{collapse}, we conclude
$
\Sigma^j C\simeq  \Sigma^j C'
$.
\end{proof}

\section{Proof of Theorem \ref{main-add}}

Suppose that $M$ and $N$ are simple homotopy equivalent $r$-connected
($r \ge 0$) closed PL manifolds of dimension $d$.

With appropriate modifications, we will
argue along the lines of Levitt's strategy for showing
$F(M,2) \simeq F(N,2)$
when $M$ and $N$ are $2$-connected (see \cite{Levitt}).

\sh{Case 1\quad $d \le 2$} 

By the classification
of low dimensional manifolds, $M$ and $N$ are PL homeomorphic. 
It follows that $F(M,k)$ and $F(N,k)$ are homeomorphic
for all $k$.

\sh{Case 2\quad $d > 2$}  
Let 
$$
\Delta_k^{\text{fat}}(M) \subset M^{\times k}
$$
denote the {\it fat diagonal.} This subpolyhedron
is the space
of $k$-tuples of points of $M$ such that at least two entries
in the $k$-tuple coincide.

By choosing a compact regular neighborhood $V \subset M^{\times k}$
of the fat diagonal, we obtain
an embedding up to homotopy of the inclusion
$\Delta_k^{\text{fat}}(M) \subset M^{\times k}$. Its complement
$C$ is weak equivalent
to $F(M,k)$ when the latter is considered as an object
of $\Top_{/M^{\times k}}$.
Denote this embedding up to homotopy
by $(V,h)$. 

Then we obtain a manifold triad  
$$(M^{\times k};V,C;\partial V)\, 
$$
(this notation means that $M^{\times k}$ is expressed
as a union of the submanifolds $V$ and $C$, with 
$V \cap C = \partial V = \partial C$).

Repeat this procedure for the fat diagonal
of $N$ in $N^{\times k}$ to get an embedding up to homotopy
of the inclusion $\Delta_k^{\text{fat}}(N) \subset N^{\times k}$.
Call the latter embedding up to homotopy $(W,h')$. Its complement
$D$ is identified with $F(N,k) \in   \Top_{/N^{\times k}}$.
Thus we have another manifold triad 
$$
(N^{\times k};W,D;\partial W)\, .
$$
The next step is to choose a 
simple homotopy equivalence $g\:N\overset\sim\to M$.
The $k$-fold product of $g$ with itself produces simple 
homotopy equivalence of pairs
$$
g_k\: (N^{\times k},\Delta_k^{\text{fat}}(N)) 
\overset\sim\to 
(M^{\times k},\Delta_k^{\text{fat}}(M))\, .
$$ 
Using the Browder-Casson-Sullivan-Wall 
theorem \cite[Th.\ 12.1]{Wall_2nd} applied to $g_k\: N^{\times k}
\to M^{\times k}$ 
and the triad  $(M^{\times k};V,C;\partial V)$,
there exists another manifold triad decomposition
of $N^{\times k}$, say
$$
(N^{\times k};V',C';\partial V') \, ,
$$
and a simple homotopy equivalence of triads 
$$
\phi\:   (N^{\times k};V',C';\partial V')
\overset\sim \to  (M^{\times k};V,C;\partial V)
$$ such that 
$\phi\:N^{\times k} \to M^{\times k}$ 
is homotopic to $g_k$.
These data describe another embedding up to homotopy of 
the inclusion $\Delta^{\text{fat}}_k(N) \to  N^{\times k}$
with the property that its complement $C'$ is identified with $F(M,k)$
up to  homotopy equivalence.

Summarizing thus far, we have two embeddings up to homotopy
of the inclusion   $\Delta^{\text{fat}}_k(N)  \to N^{\times k}$,
one whose complement is identified with $F(N,k)$ and the other
whose complement is identified with $F(M,k)$.

The next step of the argument is to verify the hypotheses
of Proposition \ref{susp-comp}. 
One checks by elementary means that 
$\dim \Delta_k^{\text{fat}}(N) \le (k{-}1)d$. 
As $d > 2$,
the hypothesis $(k{-}1)d \le kd - 3$ is satisfied.
Furthermore, the inclusion map
  $\Delta^{\text{fat}}_k(N)  \to N^{\times k}$ is $r$-connected
(recall that $r$ is the connectivity of $N$). 
Hence, applying \ref{susp-comp}
we infer 
$$
\Sigma^j D \,\, \simeq \,\, \Sigma^j C' \, ,
$$
where $j = \max(2(k{-}1)d  - kd  - r + 2,0)$. This is precisely
the case when $j = \max((k{-}2)d {-} r {+} 2,0) = \beta(k,d,r)$. 

Finally, recall that $D\simeq F(N,k)$ and 
$C'\simeq F(M,k)$. 
With respect to these identifications, we get
$$\Sigma^j F(M,k)\,\,  \simeq \,\, \Sigma^j F(N,k)\, . $$ 
This concludes the proof
of Theorem \ref{main-add}.

\section{Poincar\'e embeddings and the proof of Theorem \ref{main-thm}}

The proof of Theorem \ref{main-thm} 
will use the second author's work
on Poincar\'e embeddings from \cite{Klein_haef}
and \cite{Klein_haef2}. The material of this section
is not intended to be complete. For the foundations of 
the theory of Poincar\'e spaces, see \cite{Wall_pd} and \cite{Klein_pd}.

\sh{Motivation} 
Let $M$ be a manifold. If $K\subset M$
is a compact codimension zero submanifold, then we get a stratification of $M$ 
by submanifolds in which the codimension zero
stratum consists of $K \amalg \text{cl}(M - K)$ and the codimension one
stratum is $\partial K$.

The notion of Poincare embedding is a generalization of 
this  with the manifolds replaced by Poincar\'e spaces 
and the stratification replaced by ``stratification up to homotopy:''

\begin{defn} (cf.\ \cite{Klein_haef}).  
Let $K$ be a space which is homotopy
equivalent to a finite complex of dimension $k$, 
let $(M,\partial M)$  
a Poincar\'e duality space of dimension $m$ and let $f\: K \to M$ be a map.
A {\it (Poincar\'e) embedding}\footnote{The terminology used here differs 
slightly from 
the second author's  other papers.}
 of $f$ consists of a commutative diagram
of homotopy finite spaces $$
\SelectTips{cm}{}
\xymatrix{
A \ar[r] \ar[d] &C \ar[d]\\
K \ar[r]_{f} & M
}
$$
and a factorization $\partial M \to C \to M$ such that
\begin{itemize}
\item the diagram is homotopy cocartesian, i.e., the
map 
$$
K \times 0 \cup A \times [0,1] \cup C \times 1 \to M
$$
is a weak homotopy equivalence;
\item if $\bar K$ is the mapping cylinder of $A \to K$, then
$\bar K$ is an $m$-dimensional  Poincar\'e space with boundary $A$;
\item if $\bar C$ is the mapping cylinder of $A \amalg \partial M\to C$,
then $\bar C$ is an $m$-dimensional  
Poincar\'e space with boundary $A \amalg \partial M$;
\item there are (compatible) 
fundamental classes for $\bar K$ and $\bar C$ which glue to
give a fundamental class for $M$;
\item the map $A \to K$ is $(m{-}k{-}1)$-connected. 
\end{itemize}

The object  $C \in \Top_{\partial M \subset M}$ 
is called the {\it complement} of the (Poincar\'e) embedding.
\end{defn}

We next turn to the definition of concordance. Roughly,
a concordance can be envisioned
as a ``stratified h-cobordism'' between two
Poincar\'e embeddings. 

\begin{defn} (cf.\ \cite{Klein_haef2}) 
Assume we are  given maps $f_i\: K \to M$ for $i = 0,1$ 
which come equipped with embeddings having associated diagram
$$
\SelectTips{cm}{}
\xymatrix{
A_i \ar[r]\ar[d] & C_i \ar[d]\\
K \ar[r]_{f_i} &  M \, .
} 
$$
A {\it concordance} consists of an extension
of these data to a homotopy 
$$
F\: K\times [0,1] \to M \times [0,1]
$$ 
from 
$f_0$ to $f_1$, and a  
 commutative diagram
of homotopy finite pairs 
$$
\SelectTips{cm}{}
\xymatrix{
(A,A_0 \amalg A_1) \ar[r] \ar[d] & (W,C_0 \amalg C_1) \ar[d] \\
(K \times [0,1],K\times \partial [0,1]) \ar[r]_{F}  & 
(M \times [0,1],M \times \partial [0,1])
}
$$
together with a factorization
$$
((\partial M) \times [0,1],\partial M \times \partial [0,1])  \to 
(W ,C_0 \amalg C_1)
 \to (M \times [0,1], M \times \partial [0,1])
$$
such that 
\begin{itemize}
\item the diagram is homotopy cocartesian;
\item the inclusions $A_i \to A$ and $C_i \to W$ are homotopy 
equivalences.
\end{itemize}
\end{defn}

Relevant to the proof Theorem \ref{main-thm}  is the following
immediate consequence of the definition: a concordance 
produces a space $W$ and a commutative diagram
$$
\SelectTips{cm}{}
\xymatrix{
(\partial M) \times 0 \ar[r]^{\subset\quad }\ar[d]  
& (\partial M) \times [0,1]   \ar[d]   & \partial M \times 1 
\ar[l]_{\quad \supset} \ar[d]\\
C_0  \ar[r]^{\sim} \ar[d]  & W \ar[d] & C_1 \ar[l]_{\sim}\ar[d]\\
M\times 0  \ar[r]_{\subset\quad } & M \times [0,1]  & M \times 1 \, .
\ar[l]^{\quad \supset}
}
$$
where the horizontal arrows are homotopy equivalences.
In particular, {\it the complements of concordant embeddings are 
weak equivalent  objects of $\Top_{/M}$.}

The following is the key result used in the proof of Theorem \ref{main-thm}.
Note the loss of one dimension when compared with the manifold case.

\begin{thm}[{\cite[Cor.\ B]{Klein_haef2}}] \label{concord}  
Let $f_i\: K \to M$ for $i = 0,1$ be homotopic maps. 
Assume that $f_i$ 
come equipped with (Poincar\'e) embeddings. 
In addition, assume $k \le m{-}3$. 

If $f_0\:K \to M$ is $(2k {-}m{+}3)$-connected, then
the embeddings are concordant.
In particular, the embeddings have weak equivalent
complements.
\end{thm}

\sh{Decompression} 
Given a Poincar\'e  embedding diagram
$$
\SelectTips{cm}{}
\xymatrix{
A \ar[r]\ar[d] & C\ar[d] \\
K \ar[r]_{f} & M 
}
$$
together with factorization $\partial M \to C \to M$, apply
fiberwise suspension to obtain a new
embedding diagram 
$$
\SelectTips{cm}{}
\xymatrix{
\Sigma_K A \ar[r]\ar[d] & \Sigma_M C\ar[d]\\
K \times [0,1] \ar[r]_{f \times \text{id}} & M \times [0,1]
}
$$
together with factorization $\partial (M \times [0,1]) = \Sigma_M \partial M 
\to \Sigma_M C \to M \times [0,1]$. This operation is called {\it decompression}.

If we identify $K \times [0,1]$ with $K$ via
first factor projection, we see that decompression
increases the codimension (i.e., $m{-}k$) of the original embedding
by one. If we iterate the procedure sufficiently many times,
we eventually get into the range where Theorem \ref{concord} applies.

Hence, using Theorem \ref{concord} together with Lemma \ref{collapse},
we infer

\begin{cor}  \label{concord-conclude} Let $f_i\: K \to M$ for $i = 0,1$ be homotopic maps. 
Assume that $f_i$ 
come equipped with Poincar\'e embeddings with complements $C_i \to M$. 
In addition, assume $f_i$ is $r$-connected and $k \le m{-}3$.

Then there is a homotopy equivalence of spaces
$$
\Sigma^j C_1 \, \simeq  \,
\Sigma^j C_2 \, ,
$$
where $j = \max(2k - m - r +3,0)$.
\end{cor}

\begin{proof}[Proof of Theorem \ref{main-thm} ]
Let $h\: M \to N$ be a homotopy equivalence of $r$-connected
closed PL manifolds of dimension $d$. 
If $d \le 2$, then $M$ and $N$ are homeomorphic and
the result is trivial. From now on, assume $d > 2$.
Let $h_k\:M^{\times k} \to N^{\times k}$ be
the $k$-fold cartesian product of $h$ with itself.

Using the same notation as in the proof of Theorem \ref{main-add},
we have manifold triads $(M^{\times k};V,C;\partial V)$ and
$(N^{\times k};W,D;\partial W)$, where $V$ is a regular neighborhood
of the fat diagonal $\Delta^{\text{fat}}_k(M)$, $W$ is a regular
neighborhood of  $\Delta^{\text{fat}}_k(N)$, $C$ is identified
with the configuration space $F(M,k)$ and $D$ with the configuration
space $F(N,k)$.

Then, using the identification
$V \simeq \Delta^{\text{fat}}_k(M)$,
the first triad can be regarded as an embedding diagram
$$
\SelectTips{cm}{}
\xymatrix{
\partial V \ar[r]\ar[d]& C\ar[d]\\
 \Delta^{\text{fat}}_k(M)  \ar[r] &M^{\times k} 
}
$$
and a similar remark applies to the other triad, to give an embedding
diagram
$$\SelectTips{cm}{}
\xymatrix{
\partial W \ar[r]\ar[d] & D\ar[d] \\
 \Delta^{\text{fat}}_k(N)  \ar[r] & N^{\times k} \, .
}
$$
Applying the homotopy equivalence of pairs
$$
h_k\:(M^{\times k}, \Delta^{\text{fat}}_k(M)) \overset \sim \to 
(N^{\times k}, \Delta^{\text{fat}}_k(N))\, ,
$$
to the bottom of the first diagram, we obtain another embedding
diagram
$$
\SelectTips{cm}{}
\xymatrix{
\partial V \ar[r]\ar[d]& C\ar[d]\\
 \Delta^{\text{fat}}_k(N)  \ar[r]& N^{\times k} \, .
}
$$
Thus far, what we have achieved is two Poincar\'e embeddings of the inclusion
$\Delta^{\text{fat}}_k(N)  \to N^{\times k}$, one having 
complement $C \to M^{\times k} \to N^{\times k}$
and the other with complement $D\to N^{\times k}$.

Applying Corollary \ref{concord-conclude}, we see that
$$
\Sigma^{\alpha(k,d,r)} C \,\, \simeq \,\, \Sigma^{\alpha(k,d,r)} D \, .
$$
The proof is now concluded by referring back to  the identifications
$C \simeq F(M,k)$ and $D \simeq F(N,k)$.
\end{proof}

\Addresses\recd

\begin{thebibliography}

\bibitem[B-C-T]{BCT}%
 B\"odigheimer, C.-F., Cohen, F., Taylor, L.: 
{On the homology of configuration spaces}.  
\newblock {\it Topology}  {\bf 28},  111-123 (1989)
\MR{0991102}


\bibitem[H]{Hirschhorn}%
Hirschhorn, P.~S.: 
{Model categories and their localizations}. 
\newblock (Mathematical Surveys and Monographs, Vol. 99). 
\newblock {Amer. Math. Soc.} 2003
\MR{0991102}

\bibitem[Kl]{Klein_haef}%
Klein, J.~R.: {Poincar\'e daulity embeddings and fiberwise homotopy theory}.
\newblock {\it Topology} {\bf 38}, 597--620 (1999)
\MR{1670412}


\bibitem[Kl2]{Klein_haef2}%
Klein, J.~R.: {Poincar\'e duality embeddings and fiberwise homotopy
theory, II}.
\newblock {\it  Quart. Jour. Math. Oxford} {\bf 53} 319--335 (2002)
\MR{1930266}


\bibitem[Kl3]{Klein_pd}%
Klein, J.~R.: {Poincar\'e duality spaces}.
\newblock {Surveys on surgery theory, Vol. 1, 135--165}
\newblock {Ann.\ of Math.\ Stud.\  {145},} Princeton Univ.\ Press 2000
\MR{1747534}


\bibitem[Kr]{Kr}%
Kriz, I.:
{On the rational homotopy type of configuration spaces}.
\newblock {\it Ann. of Math.} {\bf 139}, 227--237 (1994)
\MR{1274092}


\bibitem[L]{Levitt}%
Levitt, N.: {Spaces of arcs and configuration spaces of manifolds}.  
\newblock {\it Topology} {\bf 34},  217--230 (1995)
\MR{1308497}


\bibitem[L-S]{LS}
Longoni, R., Salvatore, P.: {Configuration spaces are not homotopy invariant}. 
\arxiv{math.AT/0401075}


\bibitem[St]{Stallings}%
Stallings, J.~R.: {Embedding homotopy types into manifolds}.
\newblock 1965 unpublished paper 
(see \url{http://math.berkeley.edu/~stall} for a TeXed version)

\bibitem[Wa1]{Wall4}%
Wall, C.~T.~C.: {Classification problems in differential topology---IV.
Thickenings}.
\newblock {\it Topology\rm} {\bf 5}, 73--94 (1966)
\MR{0192509}

\bibitem[Wa2]{Wall_2nd}%
Wall, C.~T.~C.: {Surgery on Compact Manifolds}.
Second edition. Edited and with a foreword by
A. A. Ranicki. Mathematical Surveys and Monographs 69. 
Amer. Math. Soc. 1999
\MR{1687388}

\bibitem[Wa3]{Wall_pd}%
Wall, C.~T.~C.: {Poincar\'e complexes: I}.
\newblock {\it Ann.\ Math.} {\bf 86}, 213--245 (1970) 
\MR{0217791}

\end{thebibliography}
\end{document}